\documentclass{gtart_h}

\def\ifplaintex{\expandafter\ifx\csname documentclass\endcsname\relax}

\def\gtp{{\mathsurround=0pt\it $\cal G\mskip-2mu$eometry \&\ 
$\cal T\!\!$opology $\cal P\!$ublications}}  

\def\recd{{\small Received:\qua\receiveddate\ifx\reviseddate\relax
\else\qquad Revised:\qua\reviseddate\fi\par}} 

\def\lognumber#1{\def\thelognumber{#1}}
\def\volumenumber#1{\def\thevolumenumber{#1}}
\def\volumeyear#1{\def\thevolumeyear{#1}}
\def\papernumber#1{\def\thepapernumber{#1}}
\def\pagenumbers#1#2{\def\startpage{#1}\def\finishpage{#2}}
\def\published#1{\def\publishdate{#1}}

\def\received#1{\def\receiveddate{#1}}

\def\accepted#1{\def\accepteddate{#1}}

\def\asciiauthors#1{\def\theasciiauthors{#1}}
\def\asciiaddress#1{\def\theasciiaddress{#1}}
\def\asciiemail#1{\def\theasciiemail{#1}}
\def\asciiurl#1{\def\theasciiurl{#1}}
\def\coverauthors#1{\def\thecoverauthors{#1}}
\long\def\asciiabstract#1{\long\def\theasciiabstract{#1}}

\let\\\par\let\thelognumber\relax\let\thevolumenumber\relax
\let\thepapernumber\relax\let\thevolumeyear\relax\let\startpage\relax
\let\finishpage\relax\let\publishdate\relax\let\receiveddate\relax
\let\reviseddate\relax\let\accepteddate\relax\let\theasciititle\relax
\let\theasciiauthors\relax\let\theasciiaddress\relax
\let\theasciiabstract\relax

\let\thecoverauthors\relax\let\theasciiemail\relax
\let\theasciiurl\relax

\ifplaintex
\font\logobig=cmssbx10 scaled 3836
\font\logomed=cmssbx10 scaled 2557
\else
\font\logobig=cmssbx10 scaled 4200
\font\logomed=cmssbx10 scaled 2800
\fi

\long\def\makeagttitle{   
\count0=\startpage
\agt\hfill      
\hbox to 45truept{\vbox to 0pt{\vglue -13truept{\logomed A\kern -.37em{\logobig 
T}\kern -.38em G}\vss}\hss}
\break
{\small Volume \thevolumenumber\ (\thevolumeyear)
\startpage--\finishpage\nl
Published: \publishdate}

\vglue .25truein

{\parskip=0pt\leftskip 0pt plus
1fil\def\\{\par\smallskip}{\Large\bf\thetitle}\par\medskip} \vglue
0.05truein

{\parskip=0pt\leftskip 0pt plus 1fil\def\\{\par}{\sc\theauthors}
\par\medskip}%
 
\vglue 0.03truein 

{\small\leftskip 25truept\rightskip 25truept{\bf Abstract}\stdspace\theabstract

{\bf AMS Classification}\stdspace\theprimaryclass
\ifx\thesecondaryclass\relax\else; \thesecondaryclass\fi\par
{\bf Keywords}\stdspace \thekeywords\par}\vglue 7truept

}   

\ifplaintex
\font\phead=cmsl9 scaled 950
\font\pnum=cmbx10 scaled 913
\font\pfoot=cmsl9 scaled 950
\headline{\vbox to 0pt{\vskip -4.5mm\line{\small\phead\ifnum
\count0=\startpage ISSN 1472-2739 (on-line) 1472-2747 (printed)
\hfill {\pnum\folio}\else\ifodd\count0\def\\{ }%
\ifx\theshorttitle\relax\thetitle\else\theshorttitle\fi\hfill{\pnum\folio}
\else\def\\{ and }{\pnum\folio}\hfill\ifx\theshortauthors\relax\theauthors
\else\theshortauthors\fi\fi\fi}\vss}}
\footline{\vbox to 0pt{\vglue 0mm\line{\small\pfoot\ifnum\count0=\startpage
\copyright\ \gtp\hfill\else
\agt, Volume \thevolumenumber\ (\thevolumeyear)\hfill\fi}\vss}}
\else
\font=cmsl9 scaled 1050
\font\lnum=cmbx10 
\font=cmsl9 scaled 1050
\makeatletter
\def\@oddhead{{\small\lhead\ifnum\count0=\startpage ISSN 1472-2739 
(on-line) 1472-2747 (printed)\hfill {\lnum\number\count0}\else\ifodd\count0
\def\\{ }\ifx\theshorttitle\relax \thetitle \else\theshorttitle\fi\hfill
{\lnum\number\count0}\else\def\\{ and }{\lnum\number\count0}
\hfill\ifx\theshortauthors\relax 
\theauthors\else\theshortauthors\fi\fi\fi}}\def\@evenhead{\@oddhead}
\def\@oddfoot{\small\lfoot\ifnum\count0=\startpage\copyright\ \gtp\hfill\else
\agt, Volume \thevolumenumber\ (\thevolumeyear)\hfill\fi}
\def\@evenfoot{\@oddfoot}
\makeatother
\fi
\let\maketitlepage\makeagttitle
\let\makeshorttitle\maketitlepage
\let\maketitle\maketitlepage

\newwrite\gtoutfile
\long\gdef\makeheadfile{  
{\def\\{, }\def\s{ }
\immediate\openout\gtoutfile head.xxx
\immediate\write\gtoutfile{Proxy-for: \ifx\theasciiauthors\relax
\theauthors\else\theasciiauthors\fi\s<\ifx\theasciiemail\relax\theemail\else\theasciiemail\fi>}
\immediate\write\gtoutfile{\noexpand\\}
\immediate\write\gtoutfile{Authors: \ifx\theasciiauthors\relax
\theauthors\else\theasciiauthors\fi}
{\def\\{ }\immediate\write\gtoutfile{Title: \ifx\theasciititle\relax
\thetitle\else\theasciititle\fi}}
\immediate\write\gtoutfile{Subj-class: GT or SG, GR etc}
\immediate\write\gtoutfile{MSC-class: \theprimaryclass\ifx\thesecondaryclass\relax\else, \thesecondaryclass\fi}
\immediate\write\gtoutfile{Journal-ref: Algebr. Geom. Topol. \thevolumenumber\s
(\thevolumeyear) \startpage-\finishpage}
\immediate\write\gtoutfile{Comments: Published by Algebraic and
Geometric Topology at}
\immediate\write\gtoutfile{\s\s\s  http://www.maths.warwick.ac.uk/agt/AGTVol\thevolumenumber/agt-\thevolumenumber-\thepapernumber.abs.html}
\immediate\write\gtoutfile{\noexpand\\}
\immediate\write\gtoutfile{}
\ifx\theasciiabstract\relax
\immediate\write\gtoutfile{\theabstract}\else
\immediate\write\gtoutfile{\theasciiabstract}\fi
\immediate\write\gtoutfile{}
\immediate\write\gtoutfile{\noexpand\\}
\immediate\write\gtoutfile{}
\immediate\closeout\gtoutfile}}  

\def\maketitlepage{\makeagttitle\makeheadfile}
\let\makeshorttitle\maketitlepage
\let\maketitle\maketitlepage

\lognumber{33}
\volumenumber{5}
\volumeyear{2005}
\papernumber{33}
\pagenumbers{769}{784}
\received{21 January 2005} 
\accepted{31 May 2005}
\published{23 July 2005}

\usepackage{amsmath, amssymb, overpic}

\newtheorem{thm}{Theorem}[section]
\newtheorem*{claim}{Claim}

\newtheorem{prop}[thm]{Proposition}
\newtheorem{lemma}[thm]{Lemma}

\theoremstyle{remark}
\newtheorem{q}[thm]{Question}
\newtheorem{rmk}[thm]{Remark}

\newcommand{\Z}{\mathbb{Z}}

\newcommand{\bdry}{\partial}

\newcommand{\be}{\begin{enumerate}}
\newcommand{\ee}{\end{enumerate}}

\begin{document}

\title{Pinwheels and bypasses}

\shortauthors{Honda, Kazez and Mati\'c}
\authors{Ko Honda\\William H. Kazez\\Gordana Mati\'c}
\asciiauthors{Ko Honda, William H. Kazez and Gordana Matic}
\coverauthors{Ko Honda\\William H. Kazez\\Gordana Mati\noexpand\'c}

\address{{\rm KH:\qua}University of Southern California, Los Angeles, CA 90089, 
USA\\{\rm and}\\{\rm WHK and GM:\qua}University of Georgia, Athens, GA 30602, USA}

\asciiaddress{KH: University of Southern California, Los Angeles, CA 90089, 
USA\\and\\WHK and GM: University of Georgia, Athens, GA 30602, USA}

\gtemail{\mailto{khonda@math.usc.edu}, 
\mailto{will@math.uga.edu}, mailto{gordana@math.uga.edu}}

\asciiemail{khonda@math.usc.edu, will@math.uga.edu, gordana@math.uga.edu}

\gturl{\url{http://almaak.usc.edu/~khonda},
\url{http://www.math.uga.edu/~will}, \url{http://www.math.uga.edu/~gordana}} 

\asciiurl{http://almaak.usc.edu/ khonda,
http://www.math.uga.edu/ will, http://www.math.uga.edu/ gordana} 

\begin{abstract} 
We give a necessary and sufficient condition for the addition of a
collection of disjoint bypasses to a convex surface to be universally
tight -- namely the nonexistence of a polygonal region which we call a
{\em virtual pinwheel}. \end{abstract}

\asciiabstract{%
We give a necessary and sufficient condition for the addition of a
collection of disjoint bypasses to a convex surface to be universally
tight -- namely the nonexistence of a polygonal region which we call a
virtual pinwheel.}


\maketitle

\section{Introduction}

In this paper we assume that our 3-manifolds are oriented and our
contact structures cooriented. Let $\Sigma$ be a convex surface, i.e.,
it admits a $[-\varepsilon,\varepsilon]$-invariant contact
neighborhood $\Sigma\times [-\varepsilon,\varepsilon]$, where
$\Sigma=\Sigma\times \{0\}$. We do not assume that a convex surface
$\Sigma$ is closed or compact, unless specified. According to a
theorem of Giroux~\cite{Giroux00}, if $\Sigma\not=S^2$ is closed or
compact with Legendrian boundary, then $\Sigma$ has a tight
neighborhood if and only if its dividing set $\Gamma_\Sigma$ has no
homotopically trivial closed curves. (In the case when $\Sigma$ is not
necessarily compact, $\Sigma$ has a tight neighborhood if
$\Gamma_\Sigma$ has no homotopically trivial dividing curves, although
the converse is not always true.) In this paper we study the
following:

\begin{q}
Suppose we attach a family of bypasses
$\mathcal{B}=\{\mathcal{B}_\alpha\}_{\alpha \in A}$ along a disjoint
family of Legendrian arcs $\mathcal{C}=\{\delta_\alpha\}_{\alpha\in
A}$ to a product tight contact structure on $\Sigma \times
[-\varepsilon,\varepsilon]$. When is the resulting contact manifold
tight? \end{q}

\begin{figure}[ht]\small
\vspace{0.05in}
\centerline{
\begin{overpic}[height=2cm]{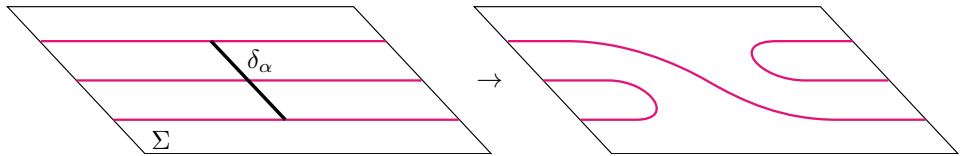}
\put(15,0.5){$\Sigma$}
\put(25,9){$\delta_\alpha$}
\put(49,7){$\to$}
\end{overpic}}
\caption{Adding a bypass}
\label{bypass}
\end{figure}

A closed Legendrian arc $\delta_\alpha$, along which a bypass
$\mathcal{B}_\alpha$ for $\Sigma$ is attached, is called a {\em
Legendrian arc of attachment}. Every arc of attachment begins and ends
on $\Gamma_\Sigma$ and has three intersection points with
$\Gamma_\Sigma$, all of which are transverse. In this paper all
bypasses are assumed to be attached ``from the front'', i.e., attached
along $\Sigma\times \{\varepsilon\}$ from the exterior of
$\Sigma\times [-\varepsilon, \varepsilon]$, and all arcs of attachment
are assumed to be embedded, i.e., there are no ``singular
bypasses''. Recall from \cite{H1} that attaching $\mathcal{B}_\alpha$
from the front and isotoping the surface $\Sigma$ across the bypass is
locally given by Figure~\ref{bypass}. Denote by $(\Sigma,
\mathcal{C})$ the contact manifold
$(\Sigma\times[-\varepsilon,\varepsilon])\cup (\cup_\alpha
N(\mathcal{B}_\alpha))$, where $N(\mathcal{B}_\alpha)$ is an invariant
neighborhood of $\mathcal{B}_\alpha$.

We will show that the key indicator of overtwistedness in the
resulting contact manifold $(\Sigma,\mathcal{C})$ is a polygonal
region in $\Sigma$ called a {\em pinwheel}. First consider an embedded
polygonal region $P$ in $\Sigma$ whose boundary consists of $2k$
consecutive sides $\gamma_1, \alpha_1, \gamma_2, \alpha_2, \dots,
\gamma_k, \alpha_k$ in counterclockwise order, where each $\gamma_i$
is a subarc of $\Gamma_\Sigma$ and each $\alpha_i$ is a subarc of a
Legendrian arc of attachment $\delta_i \in \mathcal{C}$. Here $k\geq
1$. In this paper, when we refer to a ``polygon'', we will tacitly
assume that it is an embedded polygonal region of the type just
described. Now orient the sides using the boundary orientation of
$P$. A {\em pinwheel} is a special type of polygon $P$, where, for
each $i=1,\dots,k$, $\delta_i$ extends past the final point of
$\alpha_i$ (not past the initial point) and does not reintersect
$P$. (If $k>1$, then this is equivalent to asking $\delta_i$ to extend
past $\gamma_{i+1}$, where $i$ is considered modulo $k$.)
Figure~\ref{onlyif} gives an example of a pinwheel.

It is easy to see that, if $\Sigma$ is closed or compact with
Legendrian boundary, then the addition of bypasses along all the arcs
of attachment of a pinwheel produces an overtwisted disk manifested by
a contractible curve in the resulting dividing set. Hence, the
nonexistence of pinwheels is a necessary condition for the new contact
structure to be tight. Essentially, we are asking that no closed,
homotopically trivial curves be created when some or all of the
bypasses are attached. We will prove that the nonexistence of
pinwheels is a sufficient condition as well, if $\Sigma$ is a disk
with Legendrian boundary.

\begin{thm}\label{main}
Let $\Sigma$ be a convex disk with Legendrian boundary and with a
tight neighborhood, and let $\mathcal{C}$ be a finite, disjoint
collection of bypass arcs of attachment on $\Sigma$. Denote by
$(\Sigma,\mathcal{C})$ the contact structure on $\Sigma \times I$
obtained by attaching to the product contact neighborhood of $\Sigma$
bypasses along all the arcs in $\mathcal{C}$. Then
$(\Sigma,\mathcal{C})$ is tight if and only if there are no pinwheels
in $\Sigma$.
\end{thm}

\begin{figure}[ht]\small
\centerline{
\begin{overpic}[height=3cm]{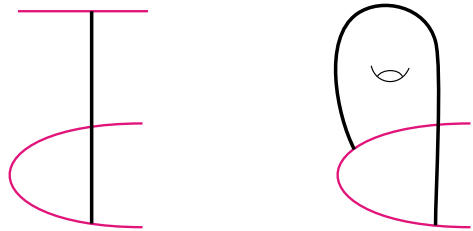} \end{overpic}}
\caption{A pinwheel and a virtual pinwheel}
\label{newtypes}
\end{figure}

If a compact convex surface $\Sigma$ has $\pi_1(\Sigma)\not=0$, then
Theorem~\ref{main} is modified to allow {\em virtual pinwheels} -- a
virtual pinwheel is an embedded polygon $P$ which becomes a pinwheel
in some finite cover of $\Sigma$. In other words, since the
fundamental group of every compact surface is residually finite, a
virtual pinwheel $P$ is either already a pinwheel or the arcs of
attachment $\delta_i$ which comprise its sides may extend beyond the
polygon, encircle a {\em nontrivial element} in $\pi_1(\Sigma,P)$ and
return to the polygon. Figure~\ref{newtypes} gives examples of arcs of
attachment which we will show result in overtwisted contact
structures. The figure to the left is a pinwheel, the one to the right
is an example of a virtual pinwheel.

\begin{thm}\label{closedsurface}
Let $\Sigma\not=S^2$ be a convex surface which is closed or compact
with Legendrian boundary and which has a tight neighborhood, and let
$\mathcal{C}$ be a finite disjoint collection of arcs of
attachment. Then the following are equivalent:
\be
\item $(\Sigma,\mathcal{C})$ is universally tight. \item There are no virtual pinwheels in $\Sigma$. \ee
\end{thm}

\begin{rmk}
A pinwheel $P$ may nontrivially intersect arcs of $\mathcal{C}$ in its
interior. Any such arc $\delta'$ would cut $P$ into two polygons, and
one of the two polygons $P'$ will satisfy the definition of a
pinwheel, with the possible exception of the condition that $\delta'$
not reintersect $P'$. $($If $\delta'$ does not reintersect $P'$, then
we can shrink $P$ to $P'$.$)$
\end{rmk}

\section{Proof of Theorem~\ref{main}}
Let $\mathcal{C}$ be the collection of arcs of attachment and
$\mathcal{B}_{\delta}$ be the bypass corresponding to
$\delta\in\mathcal{C}$.

\begin{figure}[ht]\small
\centerline{
\begin{overpic}[height=4cm]{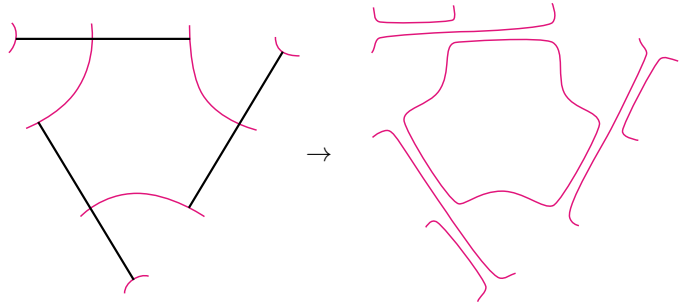}
\put(44,21){$\to$}
\end{overpic}}
\caption{Attaching the bypasses around a pinwheel}
\label{onlyif}
\end{figure}

The ``only if'' direction is immediate. If there is a pinwheel $P$,
then let $\alpha_i$, $i=1,\dots,k$, be the sides of $P$ which are
subarcs of $\delta_i\in \mathcal{C}$. Then attaching all the bypasses
$\mathcal{B}_{\delta_i}$ creates a closed homotopically trivial curve,
and hence an overtwisted disk. We can think of this disk as living at
some intermediate level $\Sigma \times \{t\}$ in the contact structure
on $\Sigma \times I$ obtained by attaching all the bypasses determined
by arcs in $\mathcal{C}$. See Figure~\ref{onlyif}.

\begin{rmk}
The arcs of attachment may be trivial arcs of attachment. For example,
if only one of them is attached, $\Gamma_\Sigma$ does not
change. However, a ``trivial arc'' $\delta_j$, after it is attached,
may affect the positions of the other arcs of attachment, if they
intersect the subarc $\gamma\subset \Gamma_\Sigma$ which forms a
polygon together with a subarc of $\delta_j$. Therefore, ``trivial''
arcs are not necessarily genuinely trivial as part of a
family. \end{rmk}

\begin{rmk}
It is crucial that an arc $\delta_i$ in the definition of a pinwheel
not return to $P$. For example, if some $\delta_j$ returns to
$\gamma_{j+1}$, then no overtwisted disk appears in a neighborhood of
the original $P$, after all the bypasses $\mathcal{B}_{\delta_i}$ are
attached.
\end{rmk}

We now prove the ``if'' part, namely if there are no pinwheels, then
the attachment of $\mathcal{C}$ onto the convex disk $\Sigma$ is
tight. In fact, we prove the following stronger result:

\begin{thm} \label{general}
Let $\Sigma$ be a convex plane whose dividing set $\Gamma_\Sigma$ has
no connected components which are closed curves. If $\mathcal{C}$ is a
locally finite, disjoint collection of bypass arcs of attachment on
$\Sigma$, and $\Sigma$ has no pinwheels, then $(\Sigma,\mathcal{C})$
is tight. \end{thm}

\begin{proof}[Reduction of Theorem~\ref{general} to Theorem~\ref{main}] 
Since any overtwisted disk\break will live in a compact region of
$\Sigma\times [-\varepsilon,\varepsilon]$, we use an exhaustion
argument to reduce to the situation where we have a closed disk $D$
with Legendrian boundary, and $\mathcal{C}$ is a finite collection of
arcs of attachment which avoid $\bdry D$. There is actually one
subtlety here when we try to use the Legendrian Realization Principle
(LeRP) on a noncompact $\Sigma$ to obtain Legendrian boundary for $D$
-- it is that there is no bound on the distance (with respect to any
complete metric on $\Sigma$) traveled by $\bdry D$ during the isotopy
given in the proof of the Giroux Flexibility Theorem. Hence we take a
different approach, namely exhausting $\Sigma$ by convex disks $D_i$
where $\bdry D_i$ is not necessarily Legendrian, and then extending
$D_i$ to a convex disk $D_i'$ with Legendrian boundary and without
pinwheels.

Let $D_1\subset D_2\subset\dots$ be such an exhaustion of $\Sigma$,
with the additional property that $\bdry D_i\pitchfork \Gamma_\Sigma$
and moreover at each intersection point $x$ the characteristic
foliation and $\bdry D_i$ agree on some small interval around
$x$. Consider a rectangle $R=[0,n]\times[0,1]$ with coordinates
$(x,y)$. Let $s$ be an arc on $\bdry D_i$ between two consecutive
intersections of $\Gamma_\Sigma\cap \bdry D_i$. Take a diffeomorphism
which takes $s$ to $x=0$; let $\xi$ be the induced contact structure
in a neighborhood of $x=0$. It is easy to extend $\xi$ to (a
neighborhood of) $y=0$ and $y=1$ so that they become dividing
curves. Now the question is to extend $\xi$ to all of $R$ so that
$x=n$ is Legendrian. Let
$R'=[0,n]\times[\varepsilon',1-\varepsilon']\subset R$ be a slightly
smaller rectangle. We write the sought-after invariant contact form on
$R'\times[-\varepsilon,\varepsilon]$ as $\alpha = dt + \beta$, where
$t$ is the coordinate for $[-\varepsilon, \varepsilon]$, $\beta$ is a
form on $R'$ which does not depend on $t$, and $d\beta$ is an area
form on $R'$. Provided $n$ is sufficiently large, $\int_{\bdry R'}
\beta$ will be positive, regardless of $\beta$ on $x=0$. Let $\omega$
be an area form on $R'$ which agrees with $d\beta$ on $\bdry R'$ and
satisfies $\int_{R'}\omega=\int_{\bdry R'}\beta$. Extend
$\beta|_{\bdry R'}$ to any 1-form $\beta'$ on $R'$ (not necessarily
the primitive of an area form). Since $d\beta'$ agrees with $\omega$
on $\bdry R'$, consider $\omega-d\beta'$. $\int_{R'}\omega-d\beta'=0$
and $\omega-d\beta'=0$ on $\bdry R'$, so by the Poincar\'e lemma there
is a 1-form $\beta''$ with $\beta''|_{\bdry R'}=0$, so that $\omega
-d\beta'=d\beta''$. Therefore, the desired $\beta$ on $R'$ is
$\beta'+\beta''$. Since there are only finitely many components of
$\Gamma_\Sigma\cap D_i$, we obtain $D'_i$ by attaching finitely many
rectangles of the type described above.
\end{proof}

Let us now consider the pair $(D,\mathcal{C})$ consisting of a convex
disk $D$ with Legendrian boundary (and dividing set $\Gamma_D$) and a
finite collection $\mathcal{C}$ of Legendrian arcs of attachment for
$D$. We now prove the following:

\begin{prop}
If $(D,\mathcal{C})$ has no pinwheels, then $(D,\mathcal{C})$ is
tight. \end{prop}

\begin{proof} The idea is to induct on the complexity of the situation. Here, the complexity $c(D,\mathcal{C})$ of $(D,\mathcal{C})$ is given by $c(D,\mathcal{C})=\#\Gamma_D+\#\mathcal{C}$, where $\#\Gamma_D$ is the number of connected components of $\Gamma_D$ and $\#\mathcal{C}$ is the number of bypass arcs in $\mathcal{C}$. Given $(D,\mathcal{C})$ we will find a pair $(D',\mathcal{C}')$ of lower complexity, where $(D,\mathcal{C})$ is tight if $(D',\mathcal{C}')$ is tight, and $(D',\mathcal{C}')$ has no pinwheels if $(D,\mathcal{C})$ has no pinwheels.

The proof will proceed by showing that there are three operations,
which we call A, B, and C, one of which can always be performed to
reduce the complexity until there are no bypasses left.  Operation A
removes (unnecessary) isolated $\bdry$-parallel arcs. If isolated
$\bdry$-parallel arcs do not exist, we apply one of Operations B and
C. If there is a trivial bypass in $\mathcal{C}$, Operation B removes
an ``innermost'' trivial bypass, i.e., we show that performing the
bypass attachment gives a configuration with lower complexity which is
tight if and only if the original configuration was tight. Otherwise,
Operation C removes an ``outermost'' nontrivial bypass by embedding
the configuration into one of lower complexity. Since each step is
strictly complexity-decreasing, and we can always do at least one of
them, we can always perform the inductive step. This will prove the
proposition and the theorem.

\subsection{Operation A: isolated $\bdry$-parallel arc}
Suppose $\Gamma_D$ has a $\bdry$-parallel arc $\gamma$ which does not
intersect any component of $\mathcal{C}$.  (Recall that arcs of
attachment are assumed to be closed and hence no component of
$\mathcal{C}$ begins or ends on $\gamma$.)  We then extend $D$ to $D'$
so that $tb(D')=tb(D)+1$ and $\Gamma_{D'}$ is obtained from $\Gamma_D$
by connecting one of the endpoints of $\gamma$ to a neighboring
endpoint of another arc in $\Gamma_D$. (See
Figure~\ref{bdryparallel}.)  It is clear that since the configuration
$(D,\mathcal{C})$ can be embedded into the configuration
$(D',\mathcal{C}'=\mathcal{C})$, and vice versa, $(D,\mathcal{C})$
tight is equivalent to $(D',\mathcal{C}')$ tight, and
$(D,\mathcal{C})$ having no pinwheels is equivalent to
$(D',\mathcal{C}')$ with no pinwheels.

\begin{figure}[ht]\small
\centerline{
\begin{overpic}[height=5cm]{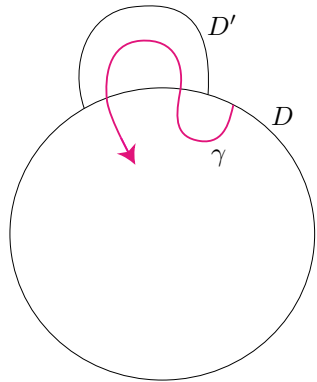}
\put(69,68){$D$}
\put(52,92){$D'$}
\put(53, 58.5){$\gamma$}
\end{overpic}}
\caption{The new $D'$}
\label{bdryparallel}
\end{figure}

\subsection{Operation B: trivial bypasses} Suppose there is a trivial arc of attachment $\delta$ in $\mathcal{C}$. Let $\gamma$ be the connected component of $\Gamma_D$ that $\delta$ intersects at least twice, and let $R$ be a closed half-disk (polygon) whose two sides are a subarc of $\delta$ and a subarc $\gamma_0$ of $\gamma$.  As shown in Figure~\ref{trivial}, we choose $R$ to be such that, with respect to the orientation of $\delta$ induced by $\bdry R$, the subarc of $\delta$ contained in $\bdry R$ starts at an interior point of $\delta$. If
$\mathcal{C}-\{\delta\}$ nontrivially intersects $int(R)$, then let
$\delta'$ be an arc of $\mathcal{C}-\{\delta\}$ which is outermost in
$R$, i.e., cuts off a subpolygon of $R$ which does not intersect
$\mathcal{C}-\{\delta\}$ in its interior.  Define $R'$ and $\gamma_0'$
analogously for $\delta'$.  (Note that $R'$ may or may not be a subset
of $R$.)  We rename $\delta'$, $R'$ and $\gamma_0'$ by omitting
primes. Therefore, we may assume that $\delta$, $\gamma$, and $R$
satisfy the property that $int(R)$ does not intersect any arc of
$\mathcal{C}$, although there may be endpoints of arcs of
$\mathcal{C}-\{\delta\}$ along $\gamma_0$. By the very definition of
$R$, the third point of intersection between $\delta$ and $\Gamma_D$
cannot be in $\gamma_0$.

\begin{figure}[ht]\small
\centerline{
\begin{overpic}[height=4.5cm]{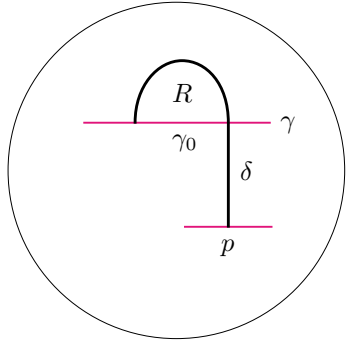}
\put(62.5,26){$p$}
\put(80,62.5){$\gamma$}
\put(48, 57){$\gamma_0$}
\put(68,47){$\delta$}
\put(48, 70){$R$}
\end{overpic}}
\caption{A trivial arc of attachment}
\label{trivial}
\end{figure}

Now, let $D=D'$ and let $\Gamma_{D'}$ be the dividing set obtained
from $\Gamma_D$ by attaching the bypass $\delta$ (the dividing set is
modified in a neighborhood of $\delta$). The isotopy type of
$\Gamma_D$ and $\Gamma_{D'}$ are the same. However, $\mathcal{C}'$ is
identical to $\mathcal{C}-\{\delta\}$ with the following exception:
arcs $\delta_i\in\mathcal{C}- \{\delta\}$ which ended on
$\gamma_0\subset \gamma$ now end on (what we may think of as) a small
interval of $\Gamma_{D'}=\Gamma_D$ around $p$. See Figures
\ref{type0shading} and \ref{type1shading}, which both depict what
happens locally near $\delta$.  We emphasize that in Figures
\ref{type0shading} and \ref{type1shading} the two dividing curves may
be part of the same dividing curve.

\begin{claim}
If $(D,\mathcal{C})$ has no pinwheels then neither does
$(D',\mathcal{C}')$. \end{claim}

\begin{proof}
For an arc $\delta_i$ in $\mathcal{C}-\{\delta\}$ with an endpoint $q$
on $\gamma_0$, let $\delta_i'$ be its image in $\mathcal{C}'$. If a
pinwheel $P'$ of $D'$ has a subarc of $\delta_i'$ as a side and $q$ as
a vertex, it is clear that there was a pinwheel $P$ of $D$ which had
subarcs of $\delta_i$ and $\delta$ as sides. The pinwheels $P$ and
$P'$ are basically the same region of $D$ -- all the sides are the
same except that $P$ has two extra vertices, $p$ and $r$, and two
extra sides. Here $r$ is the middle intersection point of $\delta$
with $\Gamma_D$ as in Figure~\ref{type0shading}.

\begin{figure}[ht]\small
\centerline{
\begin{overpic}[height=4.5cm]{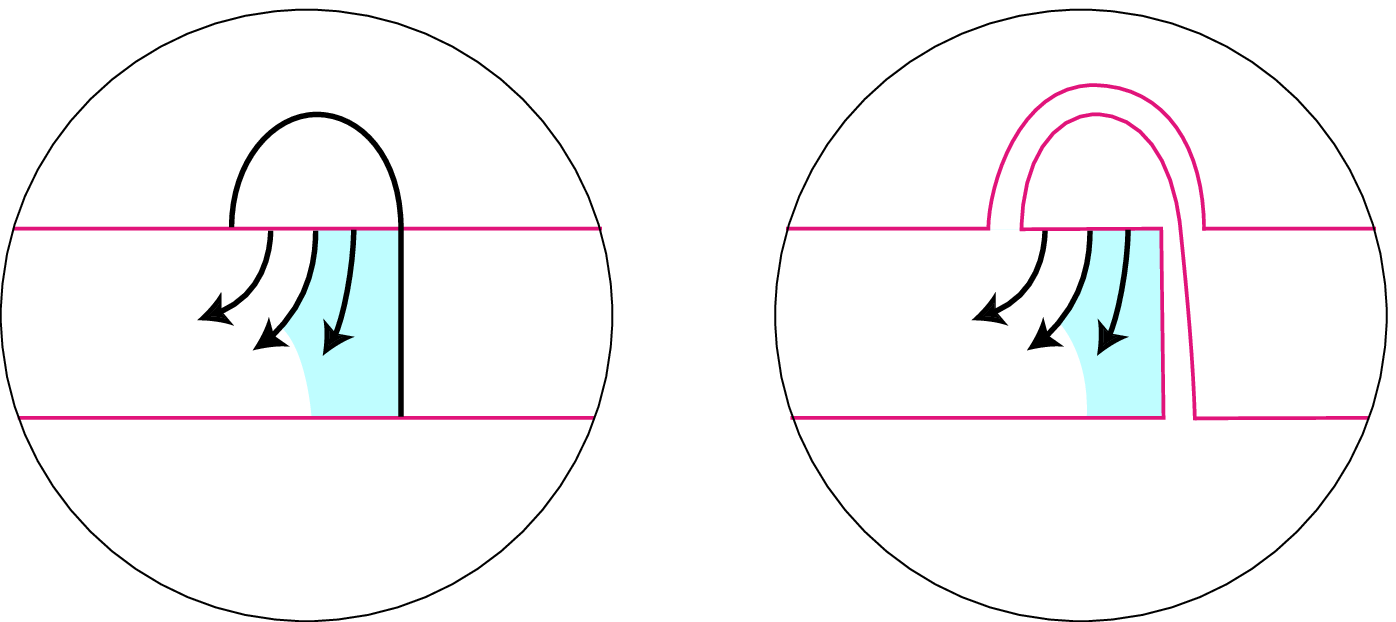}
\put(18.7,23.2){$\delta_i$}
\put(74.35,23.2){$\delta_i'$}
\put(30,12){$p$}
\put(5,30){$\gamma$}
\put(21,30){$q$}
\put(30, 30){$r$}
\put(20,37){$\delta$}
\put(24.5, 16.5){$P$}
\put(79, 16.5){$P'$}
\end{overpic}}
\caption{Pinwheels in $D$ and $D'$}
\label{type0shading}
\end{figure}

On the other hand, suppose $P'$ is a pinwheel of $D'$ which does not
involve any subarcs which used to intersect $\gamma_0$. Then $P'$ must
either completely contain or be disjoint from the region $K'$ given in
Figure~\ref{type1shading}.  However, apart from $K'$ (and the
corresponding region $K$ in $D$), $(D,\mathcal{C}-\{\delta\})$ and
$(D',\mathcal{C}')$ are identical. Hence $P'$ must have descended from
a pinwheel for $D$.
\end{proof}

\begin{figure}[ht]\small
\centerline{
\begin{overpic}[height=4.5cm]{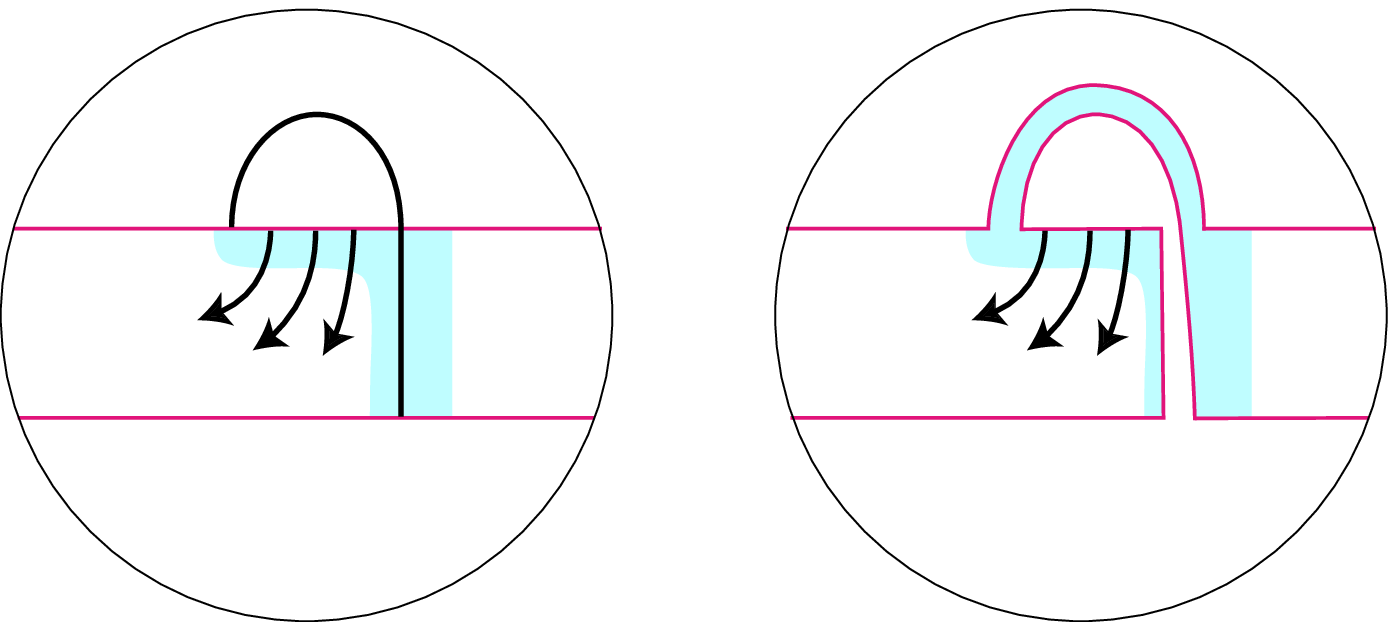}
\put(86,21){$K'$}
\put(29.25,21){$K$}
\end{overpic}}
\caption{Pinwheels in $D$ and $D'$}
\label{type1shading}
\end{figure}

\subsection{Operation C: outermost nontrivial bypass} 

Suppose all the arcs of $\mathcal{C}$ are nontrivial and there are no
isolated $\bdry$-parallel arcs. Then we have the following:

\begin{claim}
There exists an ``outermost'' arc $\delta$ with the following
property: there exists an orientation/parametrization of $\delta$ so
it intersects distinct arcs $\gamma_3$, $\gamma_2$, $\gamma_1$ of
$\Gamma_D$ in that order, and if $R\subset D$ is the closed region cut
off by the subarc $\alpha$ of $\delta$ from $\gamma_3$ to $\gamma_2$
and subarcs of $\gamma_2$ and $\gamma_3$ so that the boundary
orientation on $\alpha$ induced from $R$ and the orientation from
$\delta$ are opposite, then $int(R)$ does not intersect any other arcs
of $\mathcal{C}$. \end{claim}

\begin{figure}[ht]\small
\centerline{
\begin{overpic}[height=6cm]{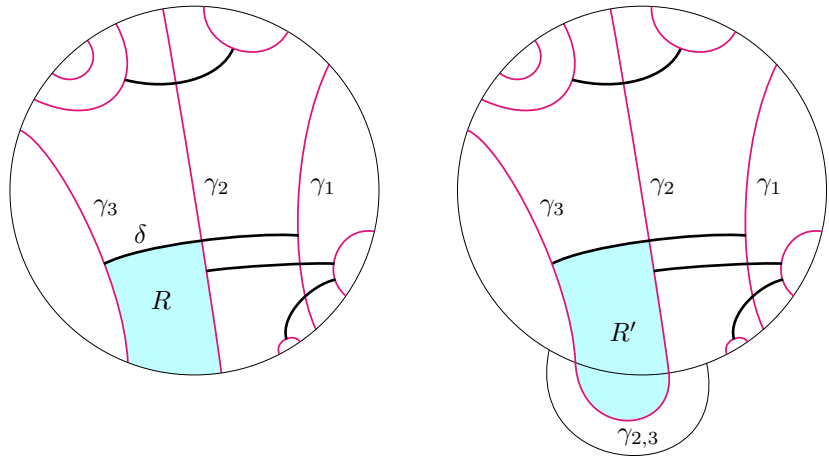}
\put(15,26){$\delta$}
\put(17,18){$R$}
\put(10, 30){$\gamma_3$}
\put(23.5, 32.5){$\gamma_2$}
\put(36.5,32.5){$\gamma_1$}
\put(64.5, 30){$\gamma_3$}
\put(78, 32.5){$\gamma_2$}
\put(91,32.5){$\gamma_1$}
\put(73,14){$R'$}
\put(74,2){$\gamma_{2,3}$}
\end{overpic}}
\caption{An outermost bypass arc}
\label{extremal}
\end{figure}

\begin{proof}
Let $\gamma$ be a $\bdry$-parallel arc of $\Gamma_D$. Since there are
no isolated $\bdry$-parallel arcs or trivial bypasses, $\gamma$
contains an endpoint of at least one arc of attachment in
$\mathcal{C}$. Of all such arcs of attachment ending on $\gamma$,
choose the ``rightmost'' one $\delta$, if we represent $D$ as the unit
disk, $\gamma$ is in the $x$-axis, and the half-disk cut off by
$\gamma$ with no other intersections with $\Gamma_D$ is in the lower
half-plane. Now orient $\delta$ so that $\gamma_3=\gamma$ and denote
by $\gamma_2$ the next arc in $\Gamma_D$ that $\delta$ intersects. Let
$R$ be the region bounded by $\gamma_3$, $\delta$, $\gamma_2$ and an
arc in the boundary of $D$, so that the boundary orientation of $R$
and the orientation of $\delta$ are opposite. There are no endpoints
of arcs of $\mathcal{C}-\{\delta\}$ along $\bdry R\cap \gamma_3$, but
there may certainly be arcs which intersect $int(R)$ and $\bdry R\cap
\gamma_2$. If there are no $\bdry$-parallel arcs in $int(R)$, then we
are done. Otherwise, take the clockwisemost $\bdry$-parallel arc
$\gamma'$ of $\Gamma_D$ (along $\bdry D$) in $int(R)$, and let
$\delta'$ be the rightmost arc of $\mathcal{C}$ starting from
$\gamma'$. Its corresponding region $R'$ is strictly contained in $R$;
hence if we rename everything by removing primes and reapply the same
procedure, then eventually we obtain $\gamma$, $\delta$, and $R$ so
that no arc of $\mathcal{C}$ intersects $int(R)$.
\end{proof}

Let $\delta$ be an outermost bypass (in the sense of the previous
claim). Then there exists an extension $D'$ of $D$, where
$tb(D')=tb(D)+1$ and $\Gamma_{D'}$ is obtained from $\Gamma_D$ by
connecting the endpoints of $\gamma_2$ and $\gamma_3$ (those which are
corners of $R$) by an arc. If $\gamma_{2,3}$ is the resulting
connected component of $\Gamma_{D'}$ which contains $\gamma_2$ and
$\gamma_3$, then $\gamma_{2,3}$ and $\delta$ cobound a disk region
$R'\subset D'$ that contains $R$. Observe that $D'$ has a tight
neighborhood, i.e., $\Gamma_{D'}$ contains no closed loops, because
$\gamma_2$ and $\gamma_3$ were distinct arcs of $\Gamma_D$. Now set
$\mathcal{C}'=\mathcal{C}$. Then $(D',\mathcal{C}')$ has lower
complexity than $(D,\mathcal{C})$. It is clear that
$(D',\mathcal{C}')$ has no pinwheels: any pinwheel $P'$ of
$(D',\mathcal{C}')$ is either already a pinwheel in $(D,\mathcal{C})$
or contains $R'$. However, any pinwheel that contains $R'$ must extend
beyond $\delta$, and hence must contain a sub-pinwheel $P$ with a side
$\delta$ inherited from $(D,\mathcal{C})$.
\end{proof}

\section{Proof of Theorem~\ref{closedsurface}} 

(1) $\Rightarrow$ (2) is clear.  Namely, if there is a virtual
pinwheel in $\Sigma$, there is a pinwheel in some finite cover, and
hence that cover is overtwisted.

(2) $\Rightarrow$ (1)\quad Assume that, on the contrary,
$(\Sigma,\mathcal{C})$ is not universally tight. Let $D$ be a disk in
the universal cover such that restriction to $D \times I$ contains the
overtwisted disk.  We can find a finite cover $\widetilde \Sigma$ of
$\Sigma$ that contains that disk, and by modifying the characteristic
foliation using LeRP if necessary, we can assume that the disk has
Legendrian boundary. Then by Theorem~\ref{main} there is a pinwheel
$P$ in $\widetilde \Sigma$. We will show that this implies the
existence of a virtual pinwheel in $\Sigma$.

\begin{figure}[ht]\small
\centerline{
\begin{overpic}[height=4.1cm]{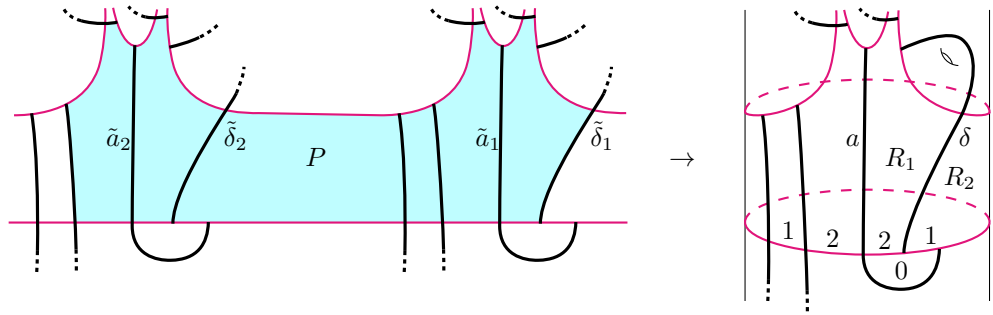}
\put(78.5,7.5){$1$}
\put(83,6.7){$2$}
\put(88.3,6.5){$2$}
\put(93,7){$1$}
\put(90,3.5){$0$}
\put(96.5,17){$\delta$}
\put(85,17){$ a$}
\put(30,15){$ P$}
\put(89,15){$R_1$}
\put(95,13){$R_2$}
\put(67,15){$\to$}
\put(21.7,17){$\tilde \delta_2$}
\put(59,17){$\tilde \delta_1$}
\put(9.7,17){$\tilde a_2$}
\put(47.3,17){$\tilde a_1$}
\end{overpic}}
\caption{Construction of a virtual pinwheel}
\label{virtualpinwheel.eps}
\end{figure}

Let $\pi: \widetilde \Sigma \rightarrow \Sigma$ be the projection map
and $\mathcal{P}$ be the set of polygons $R$ of $\Sigma$ which are
minimal in the sense that they do not contain smaller
subpolygons. Then we can define the weight function $w:
\mathcal{P}\rightarrow \Z$, which assigns to each $R\in \mathcal{P}$
the degree of $\pi^{-1}(R)\cap P$ over $R$. We illustrate this
definition in Figure~\ref{virtualpinwheel.eps}. The shaded area in the
left half of the picture is the pinwheel $P$ in $\widetilde
\Sigma$. The polygonal regions in $\Sigma$ are labeled by integers
$0,1$ and $2$ according to the value that $w$ takes on them. Our
sought-after virtual pinwheel $P'\subset \Sigma$ is then one connected
component of the union of all $R\in \mathcal{P}$ which attain the
maximal value of $w$. (In the figure there is only one component.)

First observe that $w$ cannot be locally constant, since $P$ is
strictly contained in either the positive region
$R_+(\Gamma_{\widetilde \Sigma})$ or the negative region
$R_-(\Gamma_{\widetilde \Sigma})$ -- for convenience let us suppose it
is $R_+$. Next we show that the values of $w$ are different for any
two polygonal regions $R_1$ and $R_2$ inside $R_+(\Gamma_\Sigma)$
which are adjacent along a subarc of an arc $\delta \in \mathcal{C}$
that lifts to a boundary arc $\widetilde{\delta}$ of $P$. More
precisely, orient $\delta$ so that it starts in $R_+$, and denote by
$R_1$ the region to the left of it and by $R_2$ the region to the
right (see Figure~\ref{virtualpinwheel.eps}). Then we claim that
$w(R_1) > w(R_2)$.

To prove the claim, first note that $R_1\not=R_2$. If, on the
contrary, $R_1=R_2=R$, then there is an oriented closed curve $\alpha$
such that $\alpha\backslash\delta \subset int(R)$ which is ``dual'' to
the arc $\delta$ in the sense that they intersect transversely and
$\langle \delta,\alpha\rangle = +1$. Observe that any connected
component $\widetilde {\alpha}$ of $\pi^{-1}(\alpha)\cap P$ must enter
and exit $P$ along components (say $\widetilde{\delta}_1$ and
$\widetilde{\delta}_k$) of $\pi^{-1}(\delta)$. Now, the orientation of
$\delta$ induces an orientation on the arcs in $\pi^{-1}(\delta)$, and
the induced orientation on $\widetilde{\delta}_1$ and
$\widetilde{\delta}_k$ (as seen by intersecting with
$\widetilde{\alpha}$) is inconsistent with the chirality involved in
the definition of a pinwheel.

To see what value $w$ takes on $R_1$ and $R_2$, let us look at the
components $\widetilde\delta_i$, $i=1,\dots,n$, of $\pi^{-1}(\delta)$,
and denote by $o(\delta), i(\delta)$ and $b(\delta)$ the number of
components that are respectively on the outside, in the interior, or
on the boundary of $P$. Every time a component $\widetilde\delta_i$ of
$\pi^{-1}(\delta)$ appears on the boundary of $P$, the minimal
subpolygon of $P$ adjacent to $\widetilde\delta_i$ must project to the
region $R_1$. Hence $w(R_1)= i(\delta) + b(\delta)$ and
$w(R_2)=i(\delta)$, and therefore $w(R_1) > w(R_2)$.

We will now prove that $P'$ is a polygon. We first claim that $$\pi:
\pi^{-1}(P')\cap P\rightarrow P'$$ is a covering map. Indeed, any
subpolygon of $P'$ lifts to $\max(w)$ subpolygons of $P$. Moreover, if
$a$ is a subarc of an arc $\delta'\in \mathcal{C}$, and $a$ is in
$int(P')$, then no component $\widetilde a_i$ of the lift of $a$ can
be a side of $P$, since two regions adjacent to it have equal values
of $w$. If $U$ is a neighborhood of a point on $a$ in $P'$, then
$\pi^{-1}(U) \cap P$ consists of $\max(w)$ copies of $U$. Now that we
know that $\pi^{-1}(P')\cap P$ covers $P'$, $P'$ must be simply
connected, since $\pi^{-1}(P')\cap P$ must be a union of subpolygons
of $P'$ (hence simply connected), and the cover is a finite cover.

Finally, $\delta\in \mathcal{C}$ which is a side of $P'$ and returns
to $P'$ must enclose a nontrivial element of $\pi_1(\Sigma,P')\simeq
\pi_1(\Sigma)$; otherwise it cobounds a disk together with $P'$, and
no cover of $\Sigma$ will extricate the relevant endpoint of $\delta$
from $P'$ (and hence $P$). \endproof

\begin{q}
Can we generalize Theorem~\ref{main} to the case where we have nested
bypasses? What's the analogous object to the pinwheel in this case?
\end{q}

\section{Virtual Pinwheels and Tightness} 

In this section we will discuss the following question:

\begin{q}
Can we give a necessary and sufficient condition for
$(\Sigma,\mathcal{C})$ to be tight, if $\Sigma$ is a convex surface
which is either closed or compact with Legendrian boundary?
\end{q}

\begin{figure}[ht]\small
\centerline{
\begin{overpic}[height=3cm]{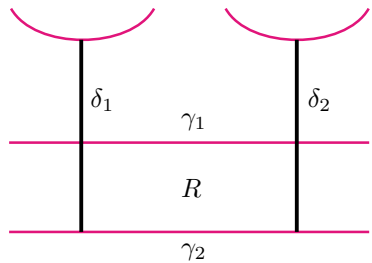}
\put(22,35){$\delta_1$}
\put(82,35){$\delta_2$}
\put(47,30){$\gamma_1$}
\put(47,-6){$\gamma_2$}
\put(47,10){$R$}
\end{overpic}}
\caption{$\delta_1$ is to the left of $\delta_2$}
\label{rotating}
\end{figure}

We will present a partial answer to this question.  Before we proceed,
we first discuss a useful technique called {\em Bypass Rotation}. Let
$\Sigma$ be a convex boundary component of a tight contact 3-manifold
$(M,\xi)$, and let $\delta_1$ and $\delta_2$ be disjoint arcs of
attachment on $\Sigma$. The bypasses are to be attached from the
exterior of $M$, and attached from the front in the figures. Suppose
there is an embedded rectangular polygon $R$, where two of the sides
are subarcs of $\delta_1$ and $\delta_2$ and the other two sides are
subarcs $\gamma_1$ and $\gamma_2$ of $\Gamma_\Sigma$. Assume
$\delta_1$ and $\delta_2$ both extend beyond $\gamma_1$ and do not
reintersect $\bdry R$. If we position $R$, $\delta_1$, and $\delta_2$
as in Figure~\ref{rotating}, so that, with the orientation induced
from $R$, $\gamma_1$ starts on $\delta_2$ and ends on $\delta_1$, then
we say that $\delta_1$ lies {\it to the left} of $\delta_2$.

The next lemma shows that the arc of attachment of a bypass can be
``rotated to the left'' and still preserve tightness. For convenience,
let $M(\delta_1,\dots,\delta_k)$ be a contact manifold obtained by
attaching $k$ disjoint bypasses from the exterior, along arcs of
attachment $\delta_1,\dots,\delta_k$.

\begin{lemma}[Bypass Rotation]\label{bypassrotating} Let $(M,\xi)$ be a contact 3-manifold, and $\delta_1, \delta_2$ be arcs of attachment on a boundary
component $\Sigma$ of $M$. If $\delta_1$ is to the left of $\delta_2$
and $M(\delta_2)$ is tight, then $M(\delta_1)$ is also
tight. \end{lemma}

\begin{proof}
If $M(\delta_2)$, is tight, then $M(\delta_1,\delta_2)=
(M(\delta_2))(\delta_1)$ is also tight, since attaching $\delta_2$
makes $\delta_1$ a trivial arc of attachment. Now,
$M(\delta_1,\delta_2)$ is also $(M(\delta_1))(\delta_2)$, so
$M(\delta_1)$ must be tight.
\end{proof}

\begin{figure}[ht!]\small
\centerline{
\begin{overpic}[height=4cm]{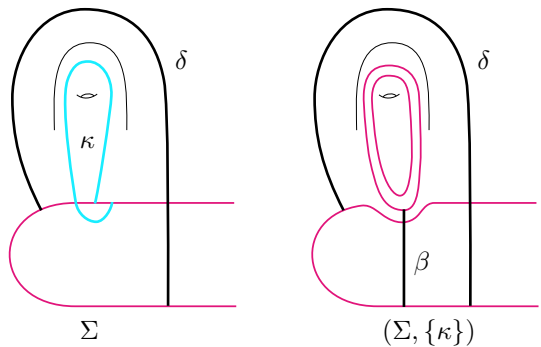}
\put(31,45){$\delta$}
\put(88,45){$\delta$}
\put(13,30){$\kappa$}
\put(76,7){$\beta$}
\put(13,-6){$\Sigma$}
\put(70,-6){$(\Sigma,\{\kappa\})$}
\end{overpic}}
\vskip.05in
\caption{Creating extra dividing curves}
\label{create}
\end{figure}

It is clear that if there is an actual pinwheel $P$ in $\Sigma$, then
$(\Sigma,\mathcal{C})$ is not tight. It is often possible to make the
same conclusion in the presence of a virtual pinwheel $P$ by using the
technique of Bypass Rotation. In fact, assume that there is a virtual
pinwheel $P$ in $\Sigma$ and let $\delta$ be an arc of attachment
which encircles a nontrivial element of $\pi_1(\Sigma, P)$. Suppose
$\delta$ has been oriented so that its orientation coincides with the
orientation on $\bdry P$.  If $\delta$ can be rotated to the left so
that the final point of $\delta$ is shifted away from $P$, then the
newly obtained configuration contains a pinwheel, and hence
$(\Sigma,\mathcal{C})$ is overtwisted by
Lemma~\ref{bypassrotating}. Even if there are no arcs of
$\Gamma_\Sigma-\bdry P$ to which we can rotate $\delta$ without
hitting other bypass arcs of attachment, we can often perform a
folding operation.  This operation can be described in two equivalent
ways (see \cite{HKM3} for details): Either add a bypass along the arc
$\kappa$ as in Figure~\ref{create} to obtain the contact structure
$(\Sigma,\mathcal{C} \cup \{\kappa\})$, or fold along a Legendrian
divide to create a pair of parallel dividing curves ``along''
$\delta$.  Since both operations can be done inside an invariant
neighborhood of $\Sigma$, $(\Sigma,\mathcal{C} \cup \{\kappa\})$ is
tight if $(\Sigma,\mathcal{C})$ is. Now, rotating $\delta$ to the
left, we can move the endpoint of $\delta$ to be on one of the newly
created parallel dividing curves -- this yields the bypass $\beta$
pictured in Figure~\ref{create}. The configuration
$(\Sigma,\mathcal{C'})$ obtained from $(\Sigma,\mathcal{C} \cup
\{\kappa\})$ by replacing $\delta$ by $\beta$ is tight if
$(\Sigma,\mathcal{C})$ is tight. However, by repeated application of
this procedure if necessary, we will often be able to eventually
obtain a genuine pinwheel $P$, hence showing that
$(\Sigma,\mathcal{C})$ is overtwisted. More precisely, we have the
following:

\begin{prop}\label{otcases}
Let $P$ be a virtual pinwheel in $\Sigma$ and $\delta\in \mathcal{C}$
be an arc of attachment on $\bdry P$ which returns to $P$.  Decompose
$\delta=\delta_0\cup\delta_1$, where $\delta_i$, $i=0,1$, have
endpoints on $\Gamma_\Sigma$ and $\delta_0\subset \bdry P$.  Orient
$\delta$ to agree with the orientation on $\bdry P$.  Let $Q$ be a
connected component of $\Sigma\setminus (\Gamma_\Sigma\cup
(\cup_{\beta\in\mathcal{C}}\beta))$ so that $\delta_1\subset \bdry Q$
and the orientation on $\delta_1$ agrees with the orientation on
$\bdry Q$.  If one of the following is true for each $\delta$, then
$(\Sigma,\mathcal{C})$ is overtwisted:
\be 
\item $Q$ is not a polygon.
\item $Q$ is a polygon but has sides in $\Gamma_\Sigma$ which are not in $\bdry P$.
\ee 
\end{prop}

\begin{figure}[ht!]\small
\centerline{
\begin{overpic}[height=5cm]{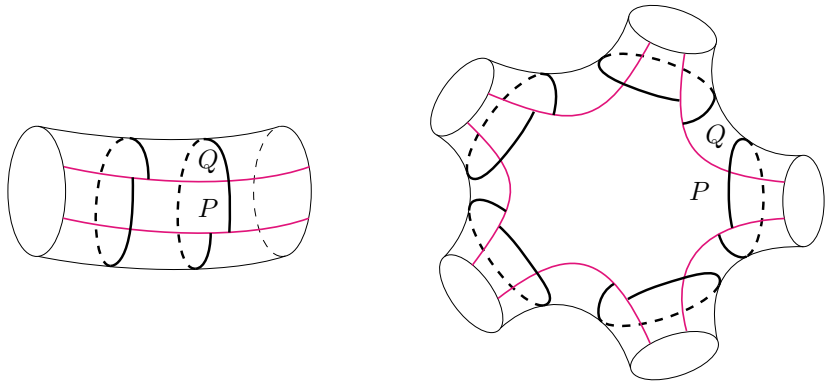}
\put(23,26){$Q$}
\put(23,20){$P$}
\put(85,29){$Q$}
\put(83,22){$P$}
\end{overpic}}
\caption{Examples of anti-pinwheels}
\label{anti-pinwheel}
\end{figure}

We now consider the situation in which the Bypass Rotation technique
just described fails.  Let $P$ be a {\it minimal pinwheel}, i.e., a
pinwheel whose interior does not intersect any arc of attachment in
$\mathcal{C}$.  Let $\delta$ be an attaching arc of $P$ that returns
to $P$ that we cannot ``unhook''. Then the region $Q$ described in
Proposition~\ref{otcases} must be polygonal and all of the edges of
$Q$ that are coming from the dividing set must be sub-edges of the
boundary of $P$.  The minimality of $P$ forces edges of $Q$ that are
attaching arcs to also be edges of $P$.  Thus $Q$ is an {\it
anti-pinwheel}, that is, a polygon whose edges which are arcs of
attachment are oriented in the direction opposite to that of a
pinwheel.  Two such pinwheel/anti-pinwheel pairs are illustrated in
Figure~\ref{anti-pinwheel}.

\begin{figure}[ht!]\small
\centerline{
\begin{overpic}[height=4cm]{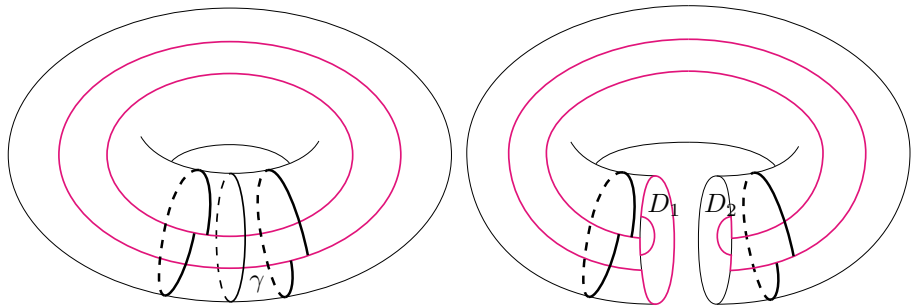}
\put(26.5,2){$\gamma$}
\put(70.5,10.3){$D_1$}
\put(76.8,10.3){$D_2$}
\end{overpic}}
\caption{Cutting the pinwheel}
\label{cut-paste}
\end{figure}

Since there is a virtual pinwheel $P$ in the situation described
above, $(\Sigma,\mathcal{C})$ is not universally tight by
Theorem~\ref{closedsurface}.  To show that there are cases when
$(\Sigma,\mathcal{C})$ is tight but virtually overtwisted, we will
analyze the situation indicated in the left portion of
Figure~\ref{cut-paste}.  Here $\Sigma = T^2$, and we think of $\Sigma
\times I$ as a neighborhood of the boundary of the solid torus. First,
we cut the solid torus along a convex disk $D$ such that its boundary
is the curve $\gamma$ that cuts $P \cup Q$ and intersects the dividing
set at two points.  Next, round the corners that are shown in
Figure~\ref{cut-paste}, and we obtain a tight convex ball. In reverse,
we can think of the solid torus as obtained by gluing two disks $D_1$
and $D_2$ on the boundary of $B^3$. The bypasses in the picture
correspond to adding trivial bypasses to the ball. By applying the
basic gluing theorem for gluing across a convex surface with
$\bdry$-parallel dividing set (see for example Theorem~1.6 in
\cite{HKM2}), we see that $(\Sigma,\mathcal{C})$ is tight.


Attempts at generalizing the above technique quickly run into some
difficulties.  Suppose we want to iteratively split $\Sigma$ along a
closed curve $\gamma$ and glue in disks $D_1$ and $D_2$.  One
difficulty (although not the only one) is that at some step in the
iteration we could get an overtwisted structure on
$(\Sigma',\mathcal{C})$.  ($\Sigma'$ could form contractible dividing
curves.)  On the other hand, this does not necessarily prove that
$(\Sigma,\mathcal{C})$ is overtwisted; it merely occurs as a subset of
a space with an overtwisted contact structure.

\rk{Acknowledgements}  KH wholeheartedly thanks the University of
Tokyo, the Tokyo Institute of Technology, and especially Prof.\
Takashi Tsuboi for their hospitality during his stay in Tokyo during
Summer-Fall 2003.  He was supported by an Alfred P.\ Sloan Fellowship
and an NSF CAREER Award. GM was supported by NSF grants DMS-0072853
and DMS-0410066 and WHK was supported by NSF grants DMS-0073029 and
DMS-0406158.  The authors also thank the referee for helpful comments.

\Addresses\recd
\end{document}